\documentclass[11pt]{article}
\usepackage[cp1251]{inputenc}
\usepackage[ukrainian, english]{babel}

\usepackage{bezier,graphics,graphicx}
\usepackage{epstopdf}
\usepackage{amsmath,amsfonts,amssymb,amscd,amsthm}
\usepackage{gensymb}
\newtheorem{theorem}{Theorem}[section]
\newtheorem{lemma}{Lemma}[section]
\newtheorem*{remark}{Remark}

\begin{document}

\title{Convergence of the directional diffusion splitting method
\footnote{
advection-diffusion-reaction model, finite element method, parabolic equation, elliptic equation, boundary value problem, recurrent integration scheme, variational problem
}
}

\author{R. Drebotiy,~ H. Shynkarenko}

\date{Ivan Franko National University of Lviv,\\
1, Universytetska St., Lviv, 79000, Ukraine
\footnotetext{UDC 518:517.948}
}
\maketitle
\thispagestyle{empty}

\begin{abstract}
We provide the proof of convergence of the directional diffusion splitting scheme for two-dimensional parabolic and elliptic advection-diffusion-reaction problems with certain restrictions on problem data.
\end{abstract}

\section{Introduction}

In the article \cite{DrebotiySplit} we introduced directional diffusion splitting scheme for two-dimensional advection-diffusion-reaction (ADR) problems of parabolic and elliptic type. This scheme allow us to parallelize computations by decomposing the each time integration step from two-dimensional problem into the set of one-dimensional ADR problems.

In this article we provide the proof of convergence of the mentioned method for parabolic problems on bounded time domain as time step goes to zero. In the same time we provide proof of convergence of the iterative process with application of the method using fixed time step as the number of steps go to infinity. In both cases there are some natural considerations on the problem data and scheme parameter choice. Also we need certain restrictions on advection field for those proofs.

The paper is structured as follows: first we define model ADR problem; then we describe the algorithm from \cite{DrebotiySplit}; after that we present mentioned proofs.

\section{Advection-diffusion-reaction problem}
Let us consider the following two-dimensional Dirichlet problem for the stationary (elliptic) advection-diffusion-reaction equation:
\begin{equation}\label{BVP1}
\left\{ \begin{aligned}
   &\text{find function}~~u:{\bar\Omega}\to \mathbb{R}~~\text{such that:}  \\
   &-\mu \Delta u+\vec{\beta }\cdot \nabla u+\sigma u=f~~~\text{in}~~\Omega\subset\mathbb{R}^2,   \\
   &u=0~~~\text{on}~~\Gamma = \partial \Omega  \\
\end{aligned} \right.
\end{equation}
and a corresponding non-stationary problem for parabolic equation:
\begin{equation}\label{Parabolic1}
\left\{ \begin{aligned}
   &\text{find function}~~u=u(x,t):{\bar\Omega \times [0,T] }\to \mathbb{R}~~\text{such that:}  \\
   &u_t^{\prime}-\mu \Delta_x u+\vec{\beta }\cdot \nabla_x u+\sigma u=f~~~\text{in}~~\Omega \times (0, T],   \\
   &u(x,t)=0,~~~(x, t)\in \Gamma \times [0, T],~~~\Gamma:=\partial\Omega \\
   &u(x,0)=u_0(x),~~~x\in \bar\Omega,  \\
\end{aligned} \right.
\end{equation}
where $\Omega $ is a bounded domain with a Lipschitz boundary $\Gamma =\partial\Omega$, $\mu =const>0$ and $\sigma =const>0$ are coefficients of diffusion and reaction respectively, function $f=f(x)$ and vector $\vec\beta =(\beta_1(x), \beta_2(x))$ represent the sources and advection flow velocity respectively. We will consider non-compressible flow, i.e., $\nabla \cdot \vec{\beta }=0$ and $\vec\beta \ne 0$ in $\Omega $. Also we consider the case, when vector field $\vec\beta$ does not have closed integral curves completely contained in $\bar\Omega$.

When we let the time $T=+\infty$ it can be shown, that the concentration $u$ stabilizes in time, i.e. we will have dynamic equilibrium, $u_t^{\prime} \rightarrow 0$ as $t\rightarrow +\infty$. In other words, we will have the degeneration of the problem \eqref{Parabolic1} to stationary problem \eqref{BVP1}.

\section{Variational formulation}

The boundary value problem \eqref{BVP1} admits the following variational formulation:
\begin{equation}\label{VP1}
\left\{ \begin{aligned}
   &\text{find}~~u\in V:=H_{0}^{1}(\Omega )~~\text{such that},  \\
   &a(u,v)=\langle l,v\rangle ~~~\forall v\in V,  \\
\end{aligned} \right.
\end{equation}
where:
\begin{equation}\label{VP1a}
\left\{ \begin{aligned}
   &a(u,v)=\int\limits_{\Omega }{(\mu \nabla_x u\cdot \nabla_x v+\vec{\beta }v\cdot \nabla_x u+\sigma uv)}dx~~~\forall u,v\in V,\\
   &\langle l,v\rangle =\int\limits_{\Omega }{f}vdx~~~\forall v\in V.\\
\end{aligned} \right.
\end{equation}
Consider standard Lebesque scalar product $(w, q) := \int_\Omega wqdx$. Using components from \eqref{VP1a} we can write the variational formulation for parabolic problem \eqref{Parabolic1}:
\begin{equation}\label{VP2}
\left\{ \begin{aligned}
   &\text{find}~~u\in L^2(0, T; V) ~~\text{such that},  \\
   &(u_t^\prime, v)+a(u,v)=\langle l,v\rangle ~~~\forall v\in V,  \\
\end{aligned} \right.
\end{equation}
Below, in this paper we will use standard $L^2$ norm of functions everywhere, unless otherwise specified. Also we consider standard Euclidean vector norm everywhere, when vector value is under the norm (if it is function, then we treat the norm of that vector-valued function as a scalar function itself).

\section{Semi-discretized directional diffusion splitting scheme}

Let us consider the normalized advection vector field $\vec b (x) = (b_1(x), b_2(x)) := \vec\beta / \|\vec\beta\|$. Consider also orthogonal vector field $\vec \gamma(x) = (\beta_2(x), -\beta_1(x))$ and corresponding normalized field $\vec p = \vec \gamma / \|\vec \gamma\|$.

Let us define time step $\Delta t$ and a fixed parameter $\theta \in (0, 1)$. Let us denote by $u_j$ (for integer $j$) an approximation to the function $u(x, t_j)$ in the time moment $t_j = j\Delta t$.

In \cite{DrebotiySplit} we propose to decompose the bilinear form of the problem into two terms
\begin{equation}\label{VPcompsplit1}
a(u,v):=s(u,v)+m(u,v),
\end{equation}
where
\begin{equation}\label{VPcompsplit2}
\left\{ \begin{aligned}
   &s(u,v):=\int\limits_{\Omega }{(\mu ({\vec b}^T\nabla_x u)({\vec b}^T\nabla_x v)+v\|\vec{\beta }\|{\vec b}^T\nabla_x u+\sigma uv)}dx\\
   &m(u,v):=\int\limits_{\Omega }{\mu ({\vec p}^T\nabla_x u)({\vec p}^T\nabla_x v)}dx\\
\end{aligned} \right.
\end{equation}
and use the following two-step scheme:
\begin{equation}\label{SemiDiscrSplit}
\left\{ \begin{aligned}
   &(\dot u_{j+\frac{1}{4}}, v) + \theta \Delta t s(\dot u_{j+\frac{1}{4}}, v) = \langle l,v\rangle - s(u_j, v) \\
   &u_{j+\frac{1}{2}} = u_j + \Delta t \dot u_{j+\frac{1}{4}},\\
   &(\dot u_{j+\frac{3}{4}}, v) + \theta \Delta t m(\dot u_{j+\frac{3}{4}}, v) = - m(u_{j+\frac{1}{2}}, v) \\
   &u_{j+1} = u_{j+\frac{1}{2}} + \Delta t \dot u_{j+\frac{3}{4}}, ~~~\forall v\in V,~~~j=0,1,...,\\
\end{aligned} \right.
\end{equation}
for obtaining approximate solution to problem \eqref{Parabolic1}. As a limit case, we also can use the procedure \eqref{SemiDiscrSplit} as an iterative algorithm for solving \eqref{VP1}. The difference is in the usage of time step $\Delta t$. For parabolic problem we need to decrease it, while keeping $T$ fixed and finite. For elliptic problems we fix $\Delta t$ and pass $j$ to infinity.

If we wish to solve original stationary problem with dominated advection, we should consider corresponding non-stationary counterpart with $u_0 \in V$ calculated as a solution of the following equation:
\begin{equation}\label{u0stationary}
s(u_0, v) = \langle l,v\rangle ~~~\forall v\in V
\end{equation}
since in that way we will obtain "from start" good approximation to the initial equation, since the dynamic of the processes will be mainly directed by advection in that case.

In \cite{DrebotiySplit} we provide parallelizable algorithm for approximate solving of each of the variational problems from the steps \eqref{SemiDiscrSplit}, and this algorithm is the key reason of possible practical usage of the scheme \eqref{SemiDiscrSplit}.

\section{Convergence study}
In this section we consider two convergence theorems: for parabolic and elliptic problems respectively. For the beginning we need one simple lemma:

\begin{lemma}\label{VectorLemma}
Consider scalar function $w: \Omega \rightarrow \mathbb{R}$: $w(x):=\|\vec\beta(x)\|$. If $\vec{\beta}\ne 0$, $\nabla \cdot \vec{\beta }=0$ and $\nabla \cdot \vec{ b }\ge 0$ in $\Omega$, then $w^\prime_{\vec\beta}\le 0$ in $\Omega$.
\end{lemma}
\begin{proof}
Consider Fig. \ref{BetaPict}.

\begin{figure}[h!]
    \centering
    \includegraphics[width=0.75\linewidth]{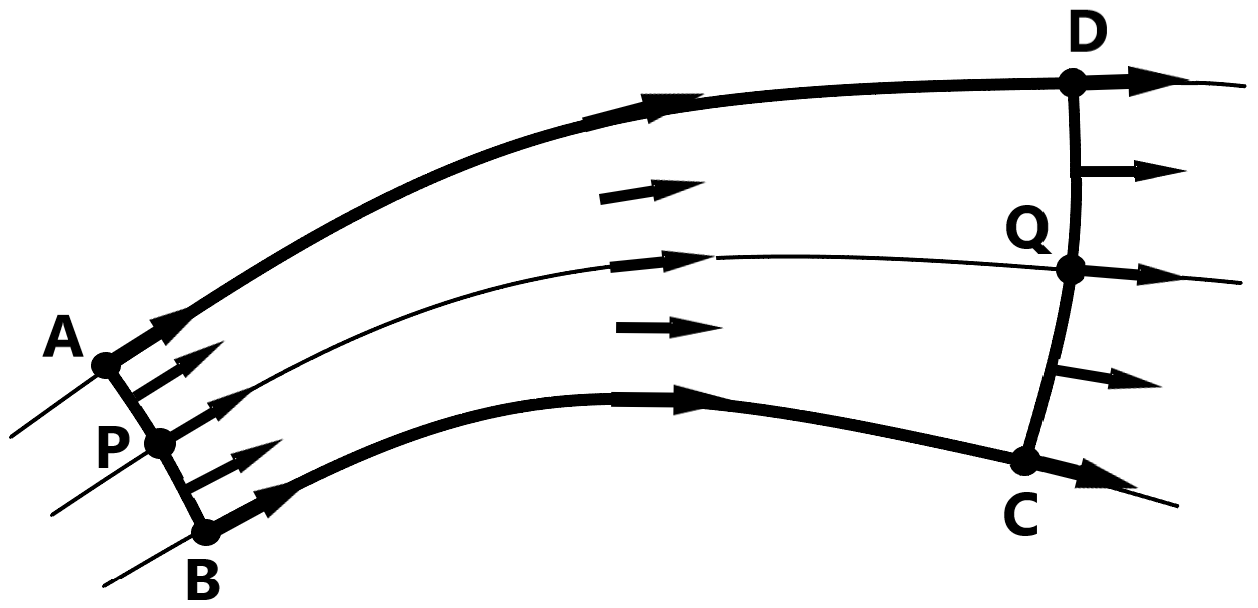}%
    \caption{Flow diagram for lemma \ref{VectorLemma}.}
    \label{BetaPict}
\end{figure}

We have there curvilinear domain $\omega:=ABCD$, spanning on the segments of flow lines of vector field $\vec\beta$ and the orthogonal field $\vec\gamma$ (segments $AB$ and $CD$). Using Gauss-Ostrogradsky and mean value theorems with the above conditions, we get:
\begin{equation}\label{lemmaEbeta}
\begin{aligned}
0&=\int\limits_{\omega}\nabla \cdot \vec{\beta }dx=\int\limits_{\partial \omega}\vec n \cdot \vec{\beta }dl=\left(\int\limits_{CD}-\int\limits_{AB}\right)\|\vec\beta\|dl=\\
&=\|\vec\beta(\xi)\| |CD| - \|\vec\beta(\eta)\| |AB|,\\
\end{aligned}
\end{equation}
where $\xi\in CD$, $\eta\in AB$ and $\vec n$ is a unit outward normal vector to the domain $\omega$.

In the same manner, taking into account that $\|\vec b\|=1$ in $\Omega$:
\begin{equation}\label{lemmaEb}
\begin{aligned}
0&\le\int\limits_{\omega}\nabla \cdot \vec{b }dx=\int\limits_{\partial \omega}\vec n \cdot \vec{b }dl=\left(\int\limits_{CD}-\int\limits_{AB}\right)dl=\\
&=|CD| - |AB|.\\
\end{aligned}
\end{equation}
By combining \eqref{lemmaEbeta} and \eqref{lemmaEb}, we obtain:
\begin{equation}\label{lemmaIneq}
\frac{\|\vec\beta(\xi)\|}{\|\vec\beta(\eta)\|} = \frac{|AB|}{|CD|} \le 1.
\end{equation}
Passing to the limit in \eqref{lemmaIneq} with $|CD|\rightarrow 0$, pulling the segments $AB$ and $CD$ to the single points $P$ and $Q$ respectively, we get
\begin{equation}\label{lemmaIneqVal}
\|\vec\beta(Q)\|\le\|\vec\beta(P)\|.
\end{equation}
Since $P$ and $Q$ are arbitrary points on the selected flow line of the vector field $\vec\beta$, we see, that the function $w(x)=\|\vec\beta(x)\|$ is not ascending along the integral curve defined by the $\vec\beta$, i.e. $w^\prime_{\vec\beta}\le 0$.
\end{proof}

\begin{theorem}
Suppose, that advection field is "divergent" in the following sense: $\nabla \cdot \vec{\beta }=0$ and $\nabla \cdot \vec{ b }\ge 0$. If the problem data ($\vec\beta$, $f$, $u_0$) and domain boundary $\partial\Omega$ are sufficiently smooth, that the classic solutions of \eqref{VP2} and variational problems from \eqref{SemiDiscrSplit} exists and parameter $\theta \ge \frac{1}{2}$, then the approximations $u_j$ defined by \eqref{SemiDiscrSplit} converge to the solution $u$ of \eqref{VP2} in $L^2$ norm as $\Delta t \rightarrow 0$. (i.e. $\Delta t \sum_{0\le j \le T/\Delta t} \|u_j-u(\cdot,t_j)\|^2 \rightarrow 0$ as $\Delta t \rightarrow 0$).
\end{theorem}
\begin{proof}
As we will see in the proof, without loss of generality, we can consider the case when $\vec\beta (x, y) = (\beta(x, y), 0)$ and $\Omega = [0,1]^2$. To transform the generic $\vec\beta$ we need to change coordinate system to other orthogonal system, when $X$ axis is going along the direction of $\vec\beta$.

From the lemma \ref{VectorLemma} we conclude, that for our setting we will have $\beta_x^\prime \le 0$.

Let us rewrite steps of \eqref{SemiDiscrSplit} by excluding intermediate values $\dot u_{j+\frac{1}{4}}$ and $\dot u_{j+\frac{3}{4}}$. We will obtain:

\begin{equation}\label{SemiDiscrSplit2}
\left\{ \begin{aligned}
   &(u_{j+\frac{1}{2}}, v) + \theta \Delta t s(u_{j+\frac{1}{2}}, v) = \Delta t \langle l,v\rangle + (u_j, v) + \Delta t (\theta - 1) s(u_j, v), \\
   &(u_{j+1}, v) + \theta \Delta t m(u_{j+1}, v) = (u_{j+\frac{1}{2}}, v) + \Delta t (\theta - 1) m(u_{j+\frac{1}{2}}, v),\\
   & ~~~\forall v\in V,~~~ 0 \le j \le T/\Delta t\\
\end{aligned} \right.
\end{equation}

Let us define operators acting on $C^2(\Omega)$:

\begin{equation}\label{SemiDiscrSplit2ClassicOperators}
\begin{aligned}
   &B_\theta:=-\theta\Delta t \mu \frac{\partial^2}{\partial x^2} + \theta\Delta t \beta (x,y) \frac{\partial}{\partial x} + (1+\theta\Delta t \sigma),\\
   &R_\theta:=-\theta\Delta t \mu \frac{\partial^2}{\partial y^2} + 1.\\
\end{aligned}
\end{equation}

In scope of our assumption on regularity, we can rewrite \eqref{SemiDiscrSplit2} in strong form:

\begin{equation}\label{SemiDiscrSplit3}
\left\{ \begin{aligned}
   &B_\theta u_{j+\frac{1}{2}}=B_{\theta-1}u_j + \Delta t f,\\
   &R_\theta u_{j+1}=R_{\theta - 1} u_{j+\frac{1}{2}}. \\
\end{aligned} \right.
\end{equation}

By excluding fractional step value $u_{j+\frac{1}{2}}$ we can rewrite \eqref{SemiDiscrSplit3} in the following way:

\begin{equation}\label{SemiDiscrSplit4}
u_{j+1}=R_\theta^{-1} R_{\theta-1} B_\theta^{-1} B_{\theta-1} u_j + \Delta t R_\theta^{-1} R_{\theta-1} B_\theta^{-1} f.
\end{equation}

To prove \textit{stability} let us estimate the norm of the operator $T:=R_\theta^{-1} R_{\theta-1} B_\theta^{-1} B_{\theta-1}$ under $L^2$ norm of the argument on the space $C^2_0(\Omega)$. We can write:

\begin{equation}\label{Tdecomp}
\begin{aligned}
   \|T\|^2\le\|R_\theta^{-1} R_{\theta-1}\|^2 \|B_\theta^{-1} B_{\theta-1}\|^2
\end{aligned}
\end{equation}
Note, that each pair of the operators $R_\theta$, $R_{\theta-1}$, $B_\theta$, $B_{\theta-1}$ commute. Let us define space $H:=\{u\in C_0^2(\Omega)|R_\theta u \in C_0^2(\Omega) \}$. Note also, that for each $v \in C_0^2(\Omega)$ we can find unique $u \in H$, such that $v=R_\theta u$. For the first norm, taking into account commutativity, we have:

\begin{equation}\label{EstimateRnorm}
\begin{aligned}
   &\|R_\theta^{-1} R_{\theta-1}\|^2 = \sup_{v\in C_0^2(\Omega)} \frac{\|R_\theta^{-1} R_{\theta-1} v\|^2}{\|v\|^2} = \sup_{u\in H} \frac{\|R_\theta^{-1} R_{\theta-1} R_\theta u\|^2}{\|R_\theta u\|^2} = \\
   &= \sup_{u\in H} \frac{\|R_\theta^{-1} R_\theta R_{\theta-1} u\|^2}{\|R_\theta u\|^2}  = \sup_{u\in H} \frac{\|R_{\theta-1} u\|^2}{\|R_\theta u\|^2} = \\
   &=\sup_{u\in H} \frac{\|-(\theta - 1)\Delta t \mu u_{yy}^{\prime\prime} + u\|^2}{\|-\theta\Delta t \mu u_{yy}^{\prime\prime} + u\|^2} = \\
   &=\sup_{u\in H} \frac{\|u\|^2 - 2 (\theta-1)\Delta t \mu (u, u_{yy}^{\prime\prime}) + (\theta-1)^2\Delta t^2\mu^2\|u_{yy}^{\prime\prime}\|^2}{\|u\|^2 - 2 \theta\Delta t \mu (u, u_{yy}^{\prime\prime}) + \theta^2\Delta t^2\mu^2\|u_{yy}^{\prime\prime}\|^2} = \\
   &=\sup_{u\in H} \frac{\|u\|^2 + 2 (\theta-1)\Delta t \mu \|u_y^\prime\|^2 + (\theta-1)^2\Delta t^2\mu^2\|u_{yy}^{\prime\prime}\|^2}{\|u\|^2 + 2 \theta\Delta t \mu \|u_y^\prime\|^2 + \theta^2\Delta t^2\mu^2\|u_{yy}^{\prime\prime}\|^2}.\\
\end{aligned}
\end{equation}

Let us assume that the last expression is less than or equal to some number $\rho\in(0,1)$. We get

\begin{equation}\label{normEstimate1}
\begin{aligned}
   &\|u\|^2 + 2 (\theta-1)\Delta t \mu \|u_y^\prime\|^2 + (\theta-1)^2\Delta t^2\mu^2\|u_{yy}^{\prime\prime}\|^2 \le \\
   &\le\rho (\|u\|^2 + 2 \theta\Delta t \mu \|u_y^\prime\|^2 + \theta^2\Delta t^2\mu^2\|u_{yy}^{\prime\prime}\|^2) \\
\end{aligned}
\end{equation}

We can regroup the last inequality in the following way:

\begin{equation}\label{normEstimate2}
\begin{aligned}
   &(1-\rho)\|u\|^2 + (\theta-1)^2\Delta t^2\mu^2\|u_{yy}^{\prime\prime}\|^2 \le \\
   &\le 2 \Delta t \mu (1-\theta(1-\rho)) \|u_y^\prime\|^2 + \rho\theta^2\Delta t^2\mu^2\|u_{yy}^{\prime\prime}\|^2 \\
\end{aligned}
\end{equation}

Taking into account Friedrichs inequality which we can write for our domain in the following form:
\begin{equation}\label{Friedrichs}
\|u\|^2 \le \frac{1}{2} \|u_y^\prime\|^2,
\end{equation}

we can see, that \eqref{normEstimate2} will be satisfied if

\begin{equation}\label{normEstimateCond}
\begin{aligned}
   &\frac{1}{2}\le\frac{2\Delta t \mu (1-\theta (1-\rho))}{1-\rho}\\
\end{aligned}
\end{equation}

and

\begin{equation}\label{normEstimateCond2}
\begin{aligned}
   &(1-\theta)^2 \le \rho \theta^2.\\
\end{aligned}
\end{equation}
Note, that \eqref{Friedrichs} and \eqref{normEstimateCond} will have different coefficient than $1/2$, if we suppose the $\Omega$ (here we denote by $\Omega$ the resulting domain after transformation of coordinate system) has more general structure and does not equal to unit square. More precisely, it will be $d^2/2$, where $d=\max\{y_2-y_1|(x,y_1),(x,y_2)\in \Omega\}$. Obviously, this will not take any effect on a further considerations.

The last inequality can be rewritten as:

\begin{equation}\label{normEstimateCond3}
\theta \ge \frac{1}{1+\sqrt{\rho}}
\end{equation}

From \eqref{normEstimateCond} and \eqref{normEstimateCond3} we see, that for any $\Delta t$ if $\theta > 1/2$, then $\| R_\theta^{-1} R_{\theta-1} \|<1$. If $\theta=1/2$, than it is not hard to see that $\|R_\theta^{-1} R_{\theta-1}\|\le 1$.

Let us estimate now $\| B_\theta^{-1} B_{\theta-1}\|$ in the similar way as in \eqref{EstimateRnorm}. Define space $G:=\{u\in C_0^2(\Omega)|B_\theta u \in C_0^2(\Omega) \}$. We have:

\begin{equation}
\begin{aligned}
   &\|B_\theta^{-1} B_{\theta-1}\|^2 = \sup_{u\in G} \frac{\|B_{\theta-1} u\|^2}{\|B_\theta u\|^2} = \\
   &=\sup_{u\in G} \frac{\|-(\theta-1)\Delta t \mu u_{xx}^{\prime\prime} + (\theta-1)\Delta t \beta u_{x}^{\prime} + (1+(\theta-1)\Delta t \sigma)u\|^2}{\|-\theta\Delta t \mu u_{xx}^{\prime\prime} + \theta\Delta t \beta u_{x}^{\prime} + (1+\theta\Delta t \sigma)u\|^2}.\\
\end{aligned}
\end{equation}

By comparing numerator and denominator, supposing the last expression is $\le 1$, with expanding the squares of the norms and using integration by parts, we can obtain the following inequality:

\begin{equation}
-(2\theta-1)\Delta t \|\mu u_{xx}^{\prime\prime}-\beta u_{x}^{\prime}\|^2-(2+(2\theta-1)\Delta t \sigma)\left[\mu\|u_{x}^{\prime}\|^2+\sigma \|u\|^2 - \frac{1}{2}(\beta_x^\prime, u^2)\right] \le 0,
\end{equation}
which is always true for $\theta \ge 1/2$ and condition $\beta_x^\prime \le 0$, showing that $\|B_\theta^{-1} B_{\theta-1} \|\le 1$. Taking into account \eqref{Tdecomp}, we have that:

\begin{equation}\label{Tnorm}
\left\{ \begin{aligned}
   &\|T\|<1, ~~~ \text{if} ~~~ \theta > 1/2,\\
   &\|T\|\le1, ~~~ \text{if} ~~~ \theta = 1/2.\\
\end{aligned} \right.
\end{equation}

Now, using \eqref{Tnorm}  and the same considerations as in \cite{Marchuk} for the representation of type \eqref{SemiDiscrSplit4}, it is not hard to see that the scheme is stable for $\theta \ge 1/2$.

To prove that the scheme is \textit{approximative} we should note, that for some starting value $u_j$ each of the variational problems on single step defined by \eqref{SemiDiscrSplit2} approximates the value of corresponding function in the time $t_{j}+\Delta t$, i.e. value of $u_{j+1}$ will be approximation to some function $w$, defined in scaled domain, i.e. in point $t_j+2\Delta t$.

Let us consider Taylor expansions for the exact solution $u(x, t)$ as a function of $t$ around the point $t=t_j$ with fixed $x$:
\begin{equation}\label{Taylor}
\left\{ \begin{aligned}
   &u\left(x, t_j+\frac{\Delta t}{2}\right)=u(x, t_j)+\frac{\Delta t}{2} u_{t}^{\prime}(x, t_j) + \frac{\Delta t^2}{8} u_{tt}^{\prime\prime}(x, \xi(x)),\\
   &u(x, t_j+\Delta t)=u(x, t_j)+\Delta t u_{t}^{\prime}(x, t_j) + \frac{\Delta t^2}{2} u_{tt}^{\prime\prime}(x, \eta(x)),\\
\end{aligned} \right.
\end{equation}

where $\xi(x) \in [t_j, t_j+\Delta t / 2]$, $\eta(x) \in [t_j, t_j+\Delta t]$. By formally substituting \eqref{Taylor} into the scheme \eqref{SemiDiscrSplit2} and considering $u$ as a function of $t$ (after substitution we will have integration over $x$) we obtain:

\begin{equation}\label{TaylorSubstitute}
\begin{aligned}
   &\frac{1}{2}(u_{t}^{\prime}(t_j), v) + s(u(t_j), v) = \langle l, v\rangle + O(\Delta t),\\
   &\frac{1}{2}(u_{t}^{\prime}(t_j), v) + m(u(t_j), v) = O(\Delta t),\\
\end{aligned}
\end{equation}

By summing up two last equations and taking into account \eqref{VPcompsplit1} we obtain:

\begin{equation}\label{TaylorSubstituteFinal}
(u_{t}^{\prime}(t_j), v) + a(u(t_j), v) = \langle l, v\rangle + O(\Delta t),
\end{equation}
as $\Delta t \rightarrow 0$. Note that we denoted by $O(\Delta t)$ some linear functional which has norm bounded with respect to $\Delta t$ as $\Delta t \rightarrow 0$

Comparing \eqref{TaylorSubstituteFinal} with exact problem \eqref{VP2} we see that the scheme \eqref{SemiDiscrSplit2} is approximative. Note, that in scope of our regularity assumptions \eqref{TaylorSubstituteFinal} can be rewritten in strong form, thus showing approximativity also in $L^2$ norm.

Now, since our scheme is approximative and stable, from the Lax equivalence theorem \cite{Lax,Marchuk,Trush} we get that it is convergent.

\end{proof}
\begin{theorem}
Suppose, that advection field is "divergent" in the following sense: $\nabla \cdot \vec{\beta }=0$ and $\nabla \cdot \vec{ b }\ge 0$. If the problem data ($\vec\beta$, $f$, $u_0$) and domain boundary $\partial\Omega$ are sufficiently smooth, that the classic solutions of \eqref{VP2} and variational problems from \eqref{SemiDiscrSplit} exists, parameter $\theta > \frac{1}{2}$ and step $\Delta t$ is fixed, then the approximation $u_j$ defined by \eqref{SemiDiscrSplit} converges to $u\in L^2(\Omega)$ as $j \rightarrow \infty$. If the limit function $u\in V$, then it is a solution of the following problem:
\begin{equation}\label{VP1limit}
\left\{ \begin{aligned}
   &\text{find}~~u\in V:=H_{0}^{1}(\Omega )~~\text{such that},  \\
   &a(u,v)=\langle l,v\rangle + \langle r, v \rangle ~~~\forall v\in V,  \\
\end{aligned} \right.
\end{equation}
where $\|r\|=O(\Delta t)$ as $\Delta t \rightarrow 0$.
\end{theorem}
\begin{proof}
Let us consider some fixed $\Delta t$. Using results from the previous proof, for $\theta > 1/2$ we have $\gamma :=\|T\|<1$. For that case for our assumptions for the process \eqref{SemiDiscrSplit4} we have $\|Tu_k - Tu_p\| \le \gamma \|u_k-u_p\|$ for $k,p=0,1,...$. Note, that in the standard proof of Banach fixed-point theorem for checking existence of the limit function, we calculate distance only between sequence elements. In that way, we can omit requirement on $T$ to be contraction on the entire definition domain, but keep it only for the set of sequence elements. Thus we will be able to proof that the sequence $\{u_j\}$ is fundamental in $L^2$ norm, and thus we have, that there exists a limit function $u\in L^2(\Omega)$. So, by passing to limit as $j\rightarrow \infty$ in \eqref{SemiDiscrSplit2}, in the case when $u\in V$, we see, that there exists also function $\tilde u \in V$, such that:
\begin{equation}\label{SemiDiscrSplitLimit}
\left\{ \begin{aligned}
   &(\tilde u, v) + \theta \Delta t s(\tilde u, v) = \Delta t \langle l,v\rangle + (u, v) + \Delta t (\theta - 1) s(u, v), \\
   &(u, v) + \theta \Delta t m(u, v) = (\tilde u, v) + \Delta t (\theta - 1) m(\tilde u, v).\\
\end{aligned} \right.
\end{equation}

First equation of the scheme can be considered as a one-$\Delta t$ time step for the initial function $u$. Taking into account convergence of such scheme with bilinear form $s(u,v)$ as $\Delta t \rightarrow 0$, we conclude that $\tilde u = u + \kappa(\Delta t)\Delta t$, where $\kappa(\Delta t)$ is bounded as $\Delta t \rightarrow 0$. By substituting this into two equations of \eqref{SemiDiscrSplitLimit} with summing them up (taking into account \eqref{VPcompsplit1}) and regrouping terms we get:

\begin{equation}\label{SemiDiscrSplitLimit2}
a(u, v)=\langle l, v \rangle + \Delta t [(\theta - 1) m(\kappa, v) - \theta s(\kappa, v)].
\end{equation}

By taking $\langle r, v \rangle := \Delta t [(\theta - 1) m(\kappa, v) - \theta s(\kappa, v)]$ we finish the proof.

\end{proof}

\begin{remark}
From the last theorem we see, that even for stationary case we should set the $\Delta t$ parameter of the iteration process sufficiently small to make final approximation matching exact stationary problem more accurate. That is why for advection-dominated problems it is naturally to choose initial approximation $u_0$ as a solution of \eqref{u0stationary}.
\end{remark}

Note, that we supposed, that the problem data is smooth enough to have existing smooth solution (which is logically to conclude after reducing problem into the set of one-dimensional problems as in \cite{DrebotiySplit}). From the nature of the processes it is reasonable also to have smooth $f$, such that $\left.f\right|_{\partial \Omega} = 0$. It is a question, in which most general class solutions $u_j$ lies. Note, that diffusion matrix in the variational forms $s(u,v)$ and $m(u,v)$ is actually projection matrix and it is in general not positive-definite, making standard proof of solution existence using Lax-Milgram lemma not possible.

Note also, that we may use a different approach to show stability. It is based on so-called energy equations \cite{Marchuk,Trush} and it can be applied to more weakened assumptions on the solutions, but, from the other hand, it imposes some limitations on scheme coefficient $\theta$. For example, let us consider fully implicit scheme, i.e. $\theta=1$. In \eqref{SemiDiscrSplit2} take $v=u_{j+\frac{1}{2}}$ in the first equation and $v=u_{j+1}$ in the second one. We will obtain:

\begin{equation}\label{SemiDiscrSplitEnergyEq}
\left\{ \begin{aligned}
   &\|u_{j+\frac{1}{2}}\|^2 + \Delta t s(u_{j+\frac{1}{2}}, u_{j+\frac{1}{2}}) = \Delta t \langle l,u_{j+\frac{1}{2}}\rangle + (u_j, u_{j+\frac{1}{2}}), \\
   &\|u_{j+1}\|^2 + \Delta t m(u_{j+1}, u_{j+1}) = (u_{j+\frac{1}{2}}, u_{j+1}).\\
\end{aligned} \right.
\end{equation}
It is not hard to see that $s(v,v)\ge \sigma \|v\|^2$ and $m(v,v)\ge 0$. Taking this into account and also Cauchy–Schwarz inequality, we obtain:
\begin{equation}\label{SemiDiscrSplitEnergyEq2}
\left\{ \begin{aligned}
   &(1+\sigma\Delta t)\|u_{j+\frac{1}{2}}\|^2 \le \Delta t \|f\| \|u_{j+\frac{1}{2}}\| + \|u_j\| \|u_{j+\frac{1}{2}}\|, \\
   &\|u_{j+1}\|^2 \le \|u_{j+\frac{1}{2}}\| \|u_{j+1}\|,\\
\end{aligned} \right.
\end{equation}
or
\begin{equation}\label{SemiDiscrSplitEnergyEq3}
\left\{ \begin{aligned}
   &(1+\sigma\Delta t)\|u_{j+\frac{1}{2}}\| \le \Delta t \|f\|  + \|u_j\| , \\
   &\|u_{j+1}\| \le \|u_{j+\frac{1}{2}}\| .\\
\end{aligned} \right.
\end{equation}

By combining inequalities from \eqref{SemiDiscrSplitEnergyEq3} we get:
\begin{equation}\label{SemiDiscrSplitEnergyEq4}
\|u_{j+1}\| \le  \frac{1}{1+\sigma\Delta t} \|u_j\| + \frac{\Delta t}{1+\sigma\Delta t} \|f\|.
\end{equation}
Since $\alpha:=1/(1+\sigma\Delta t)\le 1$ we get stability. To be more precise, we can recursively apply the \eqref{SemiDiscrSplitEnergyEq4} and obtain:
\begin{equation}\label{SemiDiscrSplitEnergyEq5}
\|u_j\|\le \alpha^j \|u_0\| + (1+\alpha + \alpha^2 + ... + \alpha^{j-1})\alpha \Delta t \|f\|.
\end{equation}
For finite $T$ (i.e. $0\le j \le T/\Delta t$) from \eqref{SemiDiscrSplitEnergyEq5} we have:
\begin{equation}\label{SemiDiscrSplitEnergyEq6}
\|u_j\|\le \|u_0\| + T \|f\|.
\end{equation}
For the case when $\sigma > 0$:
\begin{equation}\label{SemiDiscrSplitEnergyEq7}
\|u_j\|\le \alpha^j \|u_0\| + \frac{\alpha \Delta t}{1-\alpha} \|f\|.
\end{equation}
For fixed $\Delta t$ and $u_j\rightarrow u$ as $j\rightarrow \infty$ from \eqref{SemiDiscrSplitEnergyEq7} we have:
\begin{equation}\label{SemiDiscrSplitEnergyEq8}
\|u\|\le \frac{\alpha \Delta t}{1-\alpha} \|f\|.
\end{equation}

\begin{remark}
Note, that we provided proof for the splitting scheme itself, i.e. for semi-discretized scheme, in which we do not considered discretization in spatial dimensions and how it is affecting the computational process. For the fully-discretized scheme we will have iterations on finite- dimensional spaces. To study that probably we can apply known von Neumann spectral stability analysis by studying appropriate eigenvalues of matrices obtained from, for example, finite element method. In general, suppose, that we have some "projection-like" operator $D_h$, mapping exact solutions $u_j$ to appropriate discrete approximation $u_j^h$. Note, that such operator gives us result of application of spatial discretization two times, since we have two substeps on each time step. It is natural to assume, that $\|I-D_h\|\le C h$, where $h$ is some characteristic discretization parameter (like finite element diameter). In that case $\|D_h\|=1+O(h)$. If we rewrite \eqref{SemiDiscrSplit4} as $u_{j+1}=Tu_j+\Delta t L f$ and use $D$ to that scheme, we will build fully discretized variant $u_{j+1}^h=D_hTu_j^h+\Delta t D_hL f$. If showing approximativity can be done in the usual way, for stability we need to have $\|D_h T\|=\|D_h\| \|T\| \le 1$. In scope of \eqref{Tnorm}, if $\theta > 1/2 $, we need to choose $h$ small enough, to have $\|D_h\|\le 1/\|T\|$.
\end{remark}

\section{Conclusions}

In this paper under certain natural assumptions we provide the proof of convergence of the directional diffusion splitting scheme for two-dimensional parabolic and elliptic advection-diffusion-reaction problems. Several questions are open for now. First: how the error from intermediate interpolation (technically used to re-map solutions between the two substeps of the scheme in \cite{DrebotiySplit} and is not encountered in this article, as it is possible to avoid it) is propagated through the time? Second question: how to proof the convergence in the case when advection speed is variable along the flow line? Third question: how to proof convergence or apply the method in the case advection field has more complex structure or the space dimension number is greater than 2? Fourth question: how the discretization in spatial variables will affect the convergence properties (since we proved convergence of splitting procedure itself, i.e. semi-discretized approximations in time, without considering spatial discretization which can be done by different ways)?


\begin{thebibliography}{99}

\bibitem{Lax}
    \emph{Lax P.D.}
    Functional analysis // John Wiley and Sons Inc., New York, 2002.
\bibitem{Marchuk}
    \emph{Marchuk G.I., Bjorck A., Thomee V. }
    Handbook of Numerical Methods. Volume 1. Splitting and alternating direction methods. // Edited by P.G. Ciarlet and J.L. Lions. 1990, Elsevier Science Publishers B.V. (North-Holland)
\bibitem{DrebotiySplit}
    \emph{Drebotiy\,R., Shynkarenko\,H.}
    Directional diffusion splitting method for advection-diffusion-reaction model // [preprint] arXiv:2411.11041 [math.NA] 
https://doi.org/10.48550/arXiv.2411.11041
\bibitem{Trush}
    \emph{Trushevsky\,V.\,M., Shynkarenko\,H.\,A., Shcherbyna\,N.\,M.}
    Finite element method and artificial neural networks: theoretical aspects and application. //
    Lviv: Ivan Franko National University of Lviv, 2014, ISBN 978-617-10-0127-5 (in Ukrainian)
\end{thebibliography}
\end{document}